\documentclass[11pt]{article}
\usepackage{epsfig}
\usepackage{tikz}
\usetikzlibrary{cd}
\usepackage{graphicx}
\usepackage{amsthm}
\usepackage{amssymb}
\usepackage{caption}
\usepackage{subcaption}
\usepackage{color}
\usepackage{indentfirst}
\usepackage{amsmath,mathrsfs,amsfonts,verbatim,enumitem,leftidx}
\usepackage[sort&compress,numbers]{natbib}
\usepackage{amsmath,latexsym,epsf,mathabx}
\usepackage{changepage}
\usepackage[T1]{fontenc}
\usepackage[utf8]{inputenc}
\usepackage[all]{xy}
\usepackage[hidelinks]{hyperref}
\begin{document}
\newcommand{\nc}{\newcommand}
\def\PP#1#2#3{{\mathrm{Pres}}^{#1}_{#2}{#3}\setcounter{equation}{0}}
\def\ns{$n$-star}\setcounter{equation}{0}
\def\nt{$n$-tilting}\setcounter{equation}{0}
\def\Ht#1#2#3{{{\mathrm{Hom}}_{#1}({#2},{#3})}\setcounter{equation}{0}}
\def\qp#1{{${(#1)}$-quasi-projective}\setcounter{equation}{0}}
\def\mr#1{\mathrm{#1}}
\def\mc#1{{{\mathcal{#1}}}\setcounter{equation}{0}}
\def\HD{\mr{Hom}}
\def\HC{\mr{Hom}_{\mc{C}}}
\def\AdT{\mr{Add}_{\mc{T}}}
\def\adT{\mr{add}_{\mc{T}}}
\def\Kb{\mc{K}^b(\mr{Proj}R)}
\def\kb{\mc{K}^b(\mc{P}_R)}
\def\AdpC{\mr{Adp}_{\mc{C}}}
\def\AdpD{\mr{Adp}_{\mc{D}}}

\newtheorem{Th}{Theorem}[section]
\newtheorem{Def}[Th]{Definition}
\newtheorem{Lem}[Th]{Lemma}
\newtheorem{Pro}[Th]{Proposition}
\newtheorem{Cor}[Th]{Corollary}
\newtheorem{Rem}[Th]{Remark}
\newtheorem{Exm}[Th]{Example}
\newtheorem{Sc}[Th]{}
\newcommand{\cA}{\mathcal{A}}
\newcommand{\cB}{\mathcal{B}}
\newcommand{\cC}{\mathcal{C}}
\newcommand{\cF}{\mathcal{F}}
\newcommand{\cT}{\mathcal{T}}
\newcommand{\cX}{\mathcal{X}}
\newcommand{\cY}{\mathcal{Y}}
\newcommand{\sF}{\mathsf{F}}
\newcommand{\sG}{\mathsf{G}}
\newcommand{\sH}{\mathsf{H}}
\newcommand{\sT}{\mathsf{T}}
\newcommand{\Proj}{\operatorname{Proj}}
\newcommand{\add}{\operatorname{add}}
\newcommand{\pd}{\operatorname{pd}}
\newcommand{\gldim}{\operatorname{gl.dim}}
\newcommand{\id}{\operatorname{id}}
\newcommand{\Hom}{\mr{Hom}}
\newcommand{\End}{\operatorname{End}}
\newcommand{\Ext}{\mr{Ext}}
\newcommand{\Tor}{\operatorname{Tor}}
\newcommand{\Ker}{\operatorname{Ker}}
\newcommand{\Coker}{\operatorname{Coker}}
\newcommand{\Imf}{\operatorname{Im}}
\newcommand{\modcat}{\operatorname{mod}}
\newcommand{\Modcat}{\operatorname{Mod}}
\newcommand{\Findim}{\operatorname{Findim}}
\newcommand{\findim}{\operatorname{findim}}
\newcommand{\pamp}{\operatorname{pamp}}
\newcommand{\st}{\mathrm{st}}
\newcommand{\op}{\mathrm{op}}

\def\Pf#1{{\noindent\bf Proof}.\setcounter{equation}{0}}
\def\>#1{{ $\Rightarrow$ }\setcounter{equation}{0}}
\def\<>#1{{ $\Leftrightarrow$ }\setcounter{equation}{0}}
\def\bskip#1{{ \vskip 20pt }\setcounter{equation}{0}}
\def\sskip#1{{ \vskip 5pt }\setcounter{equation}{0}}
\def\mskip#1{{ \vskip 10pt }\setcounter{equation}{0}}

\def\bg#1{\begin{#1}\setcounter{equation}{0}}
\def\ed#1{\end{#1}\setcounter{equation}{0}}
\def\KET{T^{^F\bot}\setcounter{equation}{0}}
\def\KEC{C^{\bot}\setcounter{equation}{0}}

\def\jze{{ \begin{pmatrix} 0 & 0 \\ 1 & 0 \end{pmatrix}}\setcounter{equation}{0}}
\def\hjz#1#2{{ \begin{pmatrix} {#1} & {#2} \end{pmatrix}}\setcounter{equation}{0}}
\def\ljz#1#2{{  \begin{pmatrix} {#1} \\ {#2} \end{pmatrix}}\setcounter{equation}{0}}
\def\jz#1#2#3#4{{  \begin{pmatrix} {#1} & {#2} \\ {#3} & {#4} \end{pmatrix}}\setcounter{equation}{0}}


\title{\bf  Auslander bounds under recollements and arrow removal}
\smallskip

\author{\small Kaili Wu\\ \small College of Science, Nanjing Forestry University, Nanjing, \rm210037, P. R. China\\
Email address: kailywu@163.com\\  \small Jiaqun Wei\\ \small School of Mathematical  Sciences, Zhejiang  Normal University, Jinhua \rm321004, P. R. China\\Email address: weijiaqun5479@zjnu.edu.cn
}
\date{}
\maketitle

\baselineskip 15pt

\noindent {\bf Abstract}:  Auslander bounds measure uniformity in the eventual vanishing
of extension groups and are closely related to homological
conjectures, including the Auslander-Reiten conjecture and   finitistic dimension conjecture.
We study global and finitistic Auslander bounds of
finite-dimensional algebras under recollements and arrow
removal. For a recollement induced by an idempotent $e$
of an algebra $\Lambda$, we establish inequalities comparing
the bounds of $\Lambda$ with those of the quotient algebra
$\Lambda/\Lambda e\Lambda$ and the corner algebra $e\Lambda e$. 
We also establish comparison results for cleft extensions
and deduce that arrow removal preserves finiteness of both
the global and finitistic Auslander bounds. Finally, we construct radical-square-zero algebras with $\findim(\Lambda_n)=0$ and $\mr{fab}(\Lambda_n)=\mr{gb}(\Lambda_n)=n$ for every $n\geq1$, yielding arbitrarily large gaps between finitistic Auslander bounds and finitistic dimensions.

\mskip\


\smallskip


\section{Introduction}
Auslander's conjecture on eventual Ext vanishing asserts that,
for every Artin algebra $R$ and every finitely generated left
$R$-module $M$, there is a nonnegative integer $b_M$ such that
\[
\mr{Ext}_R^i(M,N)=0
\qquad\text{for all }i>b_M
\]
whenever $N\in R\text{-}\mr{mod}$ satisfies
$\mr{Ext}_R^j(M,N)=0$ for all sufficiently large $j$;
see \cite{H}. The bound $b_M$ may depend on $M$, but is
independent of $N$. Although the conjecture fails in general
\cite{JS,Sm}, it motivates the study of quantitative bounds
on eventual Ext vanishing. The least such nonnegative integer
is called the Auslander bound of $M$, denoted by $b_R(M)$;
we set $b_R(M)=\infty$ if no such integer exists.
An Artin algebra for which this vanishing property holds is said
to satisfy the Auslander Condition; equivalently, every finitely
generated left module has finite Auslander bound.
Rings satisfying this condition have been studied in \cite{CH,HJ,Mori}.

Let $\Lambda$ be a finite-dimensional algebra over an
algebraically closed field. The global and finitistic Auslander
bounds are defined, respectively, by
\[
\begin{aligned}
\mr{gb}(\Lambda)
&=\sup\{b_\Lambda(M)\mid M\in\Lambda\text{-}\mr{mod}\},\\
\mr{fab}(\Lambda)
&=\sup\{b_\Lambda(M)\mid M\in\Lambda\text{-}\mr{mod},\
b_\Lambda(M)<\infty\}.
\end{aligned}
\]
Thus $\mr{fab}(\Lambda)$ ranges only over modules with finite
Auslander bound, whereas $\mr{gb}(\Lambda)$ ranges over all modules.
The Uniform Auslander Condition requires
$\mr{gb}(\Lambda)<\infty$. If $\Lambda$ satisfies the
Auslander Condition, then $\mr{fab}(\Lambda)=\mr{gb}(\Lambda)$.

Since $b_\Lambda(M)=\mr{pd}_\Lambda M$ whenever
$\mr{pd}_\Lambda M<\infty$, one has
\[
\mr{findim}(\Lambda)
\leq\mr{fab}(\Lambda)
\leq\mr{gb}(\Lambda).
\]
where $\mr{findim}(\Lambda)=
\sup\{\mr{pd}_\Lambda M
\mid M\in\Lambda\text{-}\mr{mod},\
\mr{pd}_\Lambda M<\infty\}$
denotes the finitistic dimension of $\Lambda$.
In \cite{W}, the second author formulated the finitistic Auslander
conjecture, asserting that $\mr{fab}(\Lambda)<\infty$ for every
Artin algebra $\Lambda$. By the displayed inequality, this
conjecture implies the finitistic dimension conjecture.
The Auslander Condition also implies the Auslander--Reiten conjecture,
which asserts that a finitely generated left $\Lambda$-module $M$
is projective whenever $\mr{Ext}_\Lambda^i(M,M\oplus\Lambda)=0$
for all $i>0$ \citep[Theorem A]{CH}.

The behavior of Auslander bounds under equivalences provides
a starting point for comparison results. Using syzygy complexes,
the second author obtained derived-category characterizations of the finiteness
of the global and finitistic Auslander bounds and proved that
derived equivalences preserve both finiteness properties,
as well as the Auslander Condition
\citep[Theorems 6.2 and 6.3]{Wei2018}.
Auslander bounds for complexes have also been studied in \cite{L}.
These results motivate comparisons of Auslander bounds under
homological constructions beyond equivalences, particularly
recollements and arrow removal.

Recollements, introduced by Beilinson, Bernstein and Deligne
\cite{BBD}, relate three categories through a system of adjoint
functors. In the representation theory of algebras, they have
been used to compare global dimension \cite{K}, finitistic
dimension \cite{CX,PH}, and Gorensteinness \cite{QH}.
Psaroudakis and Vit\'oria \cite{PV} proved that a recollement
whose terms are module categories is equivalent to one induced
by an idempotent element.

In this paper, we study the global and finitistic Auslander
bounds of algebras related by idempotent-induced recollements.
For an idempotent $e\in\Lambda$, we consider the recollement
involving
\[
(\Lambda/\Lambda e\Lambda)\text{-}\mr{mod},
\qquad
\Lambda\text{-}\mr{mod},
\qquad
e\Lambda e\text{-}\mr{mod}.
\]
We establish explicit comparisons between the Auslander bounds
of $\Lambda$ and those of its quotient and corner algebras.
The following theorem summarizes the resulting finiteness
criteria; the explicit estimates are given in
Theorems \ref{T1} and \ref{T2}.

\begin{Th}\label{1}
Let $e$ be an idempotent element of $\Lambda$.
Then the following statements hold.
\begin{enumerate}
\item[(a)]
Suppose that $\Lambda e\Lambda$ is a strong idempotent ideal and
$
\mr{gldim}(e\Lambda e)<\infty,$
$\mr{pd}_{\Lambda}(\Lambda e\Lambda)<\infty,$ $\mr{pd}_{\Lambda^{op}}(\Lambda e\Lambda)<\infty.$
Then
\begin{align*}
\mr{gb}(\Lambda)<\infty
&\quad\Longleftrightarrow\quad
\mr{gb}(\Lambda/\Lambda e\Lambda)<\infty,\\
\mr{fab}(\Lambda)<\infty
&\quad\Longleftrightarrow\quad
\mr{fab}(\Lambda/\Lambda e\Lambda)<\infty.
\end{align*}

\item[(b)]
Suppose that $\mathbb{P}_1=\mathbb{P}_\infty$ and
$
\mr{gldim}(\Lambda/\Lambda e\Lambda)<\infty,$
$\mr{pd}_{\Lambda}(\Lambda e\Lambda)<\infty.
$
Then
\begin{align*}
\mr{gb}(\Lambda)<\infty
&\quad\Longleftrightarrow\quad
\mr{gb}(e\Lambda e)<\infty,\\
\mr{fab}(\Lambda)<\infty
&\quad\Longleftrightarrow\quad
\mr{fab}(e\Lambda e)<\infty.
\end{align*}
\end{enumerate}
\end{Th}

Here $\mathbb{P}_1$ and $\mathbb{P}_\infty$ are the
subcategories associated with $e$ defined in Section 2.
The two parts reduce the finiteness questions to the quotient
and corner algebras, respectively.
We illustrate these results with examples and obtain sharper
estimates for upper triangular matrix algebras; see
Proposition \ref{P-upper-triangular-sharp}.

We also study the arrow-removal operation introduced by
Green, Psaroudakis and Solberg \cite{GPS2021}.
Using comparison results for cleft extensions, we prove that
arrow removal preserves finiteness of both Auslander bounds.
The following theorem summarizes this result; the quantitative
estimates are given in Theorem \ref{T-arrow-bound}.

\begin{Th}\label{T-intro-arrow}
Let $\Lambda=\Bbbk Q/I$, where $Q$ is a finite quiver and
$I$ is an admissible ideal. Suppose that
$a_i:v_{e_i}\longrightarrow v_{f_i}$ $(1\leq i\leq s)$
are arrows that do not occur in a fixed minimal generating set
of $I$, and that
\[
\mr{Hom}_{\Lambda}(\Lambda f_i,\Lambda e_j)=0
\qquad(1\leq i,j\leq s).
\]
Set $\overline a_i=a_i+I$ and
\[
\Gamma=\Lambda\big/\sum_{i=1}^{s}\Lambda\overline a_i\Lambda.
\]
Then
\begin{align*}
\mr{gb}(\Lambda)<\infty
&\quad\Longleftrightarrow\quad
\mr{gb}(\Gamma)<\infty,\\
\mr{fab}(\Lambda)<\infty
&\quad\Longleftrightarrow\quad
\mr{fab}(\Gamma)<\infty.
\end{align*}
\end{Th}

The paper is organized as follows.
Section 2 recalls recollements, idempotent ideals, and the
definitions of Auslander bounds.
Section 3 proves the reduction theorems and derives sharper
estimates for upper triangular matrix algebras.
Section 4 establishes comparison results for cleft extensions
and applies them to arrow removal.
Section 5 constructs a radical-square-zero family for which
$\mr{fab}(\Lambda)-\findim(\Lambda)$ takes arbitrarily large finite values.

\section{Preliminaries}
%
%

Unless otherwise specified, $\Lambda$ is a finite-dimensional algebra
over an algebraically closed field $\Bbbk$, and modules are finitely
generated left modules. We write $\Lambda\text{-}\mr{mod}$ for the
category of finitely generated left $\Lambda$-modules. All Auslander bounds and
finitistic dimensions are taken on the left. We put
$\mathbb N=\{0,1,2,\ldots\}$ and write
$D=\mr{Hom}_{\Bbbk}(-,\Bbbk)$ for the standard duality between
finitely generated left and right $\Lambda$-modules.

For $X\in\Lambda\text{-}\mr{mod}$, we denote its projective and
injective dimensions by $\mr{pd}_\Lambda X$ and $\mr{id}_\Lambda X$,
respectively. We write $\Omega^n_\Lambda X$ and $\Omega^{-n}_\Lambda X$
for its $n$th syzygy and cosyzygy, taken from minimal projective
resolutions and minimal injective coresolutions unless otherwise
specified, and set $\Omega^0_\Lambda X=X$.
The full subcategory $\mr{add}_\Lambda X$ consists of all direct
summands of finite direct sums of copies of $X$.
For a bimodule ${}_A U_B$, the notation $\mr{pd}_A U$ and
$\mr{pd}_{B^{op}}U$ refers to its left $A$-module and right
$B$-module projective dimensions, respectively.

\subsection{Recollements and idempotent ideals}

We recall the definition and basic properties of recollements
of abelian categories.
\begin{Def}\citep[Definition 2.1]{PH}
A recollement of abelian categories $\mathcal{A}$, $\mathcal{B}$ and $\mathcal{C}$ is a diagram
\begin{eqnarray*}
\xymatrix{{\mathcal{B}}\ar[rr]|{i_*=i_!}&&\mathcal{A}\ar[rr]|{j^!=j^*}\ar[ll]<3ex>|{i^!}\ar[ll]<-3ex>|{i^*}&&\mathcal{C}\ar[ll]<3ex>|{j_*}\ar[ll]<-3ex>|{j_!}}
\end{eqnarray*}
henceforth denoted by $\mr{R}(\mathcal{A},\mathcal{B}, \mathcal{C})$, satisfying the following conditions:\\
(RA1) $(i^{*},i_{*},i^{!})$ and $(j_{!}, j^{*},j_{*})$ are  adjoint triples.\\
(RA2) The functors $i_{*}$, $j_{!}$ and $j_{*}$ are fully faithful.\\
(RA3) $\mr{Im}i_{*}= \mr{Ker}j^{*}$.
\end{Def}
A Serre subcategory of an abelian category is a full subcategory
closed under subobjects, quotients, and extensions.
The following proposition collects standard consequences of the definition.
\begin{Pro}\cite{FP, PH, PV}\label{P2.2}
The following hold for a recollement $\mr{R}(\mathcal{A},\mathcal{B}, \mathcal{C})$.\\
(a) The functors $i_{*}$ and $j^{*}$ are exact.\\
(b) The compositions $i^{*} j_{!}$ and $i^{!} j_{*}$ are identically zero.\\
(c) The units  $\mr{Id}_{\mathcal{B}}\to i^{!}i_{*}$ and $\mr{Id}_{\mathcal{C}}\to j^{*}j_{!}$ and the counits $i^{*}i_{*}\to \mr{Id}_{\mathcal{B}}$ and $ j^{*} j_{*}\to\mr{Id}_{\mathcal{C}}$ are natural isomorphisms.\\
(d) The image of $i_*$ is a Serre subcategory of $\mathcal A$, and
$j^*$ induces an equivalence
$\mathcal A/\mr{Im}i_*\simeq\mathcal C$.\\
(e) For every $A\in\mathcal A$, there exist $B,B'\in\mathcal B$
such that the units and counits of the adjunctions induce exact sequences
$$0\to i_{*}(B)\to  j_{!}j^{*}(A)\to A\to i_{*}i^{*}A\to 0,$$
$$0\to i_{*}i^{!}(A)\to A\to j_{*}j^{*}A\to i_{*}(B')\to 0\text{. }$$
\end{Pro}

\begin{Exm}
Let $e$ be an idempotent element of $\Lambda$, and put
$I=\Lambda e\Lambda$.
Since $\Lambda$ is finite-dimensional, all six functors in the
idempotent-induced recollement of categories of all modules
\cite{KU,PH} preserve finite generation. Hence this recollement
restricts to the following recollement of categories of finitely
generated left modules:
\begin{eqnarray*}
\xymatrix{{\Lambda/I\text{-}\mr{mod}}\ar[rr]^{inc}&&\Lambda\text{-}\mr{mod}\ar[rr]^{e\Lambda\otimes_\Lambda-}\ar[ll]<3ex>^{\mr{Hom}_{\Lambda}(\Lambda/I,-)}\ar[ll]<-3ex>_{\Lambda/I\otimes_\Lambda-}&&e\Lambda e\text{-}\mr{mod}\ar[ll]<3ex>^{\mr{Hom}_{e\Lambda e}(e\Lambda,-)}\ar[ll]<-3ex>_{\Lambda e\otimes_{e\Lambda e}-}}.
\end{eqnarray*}
\end{Exm}

Strongly idempotent ideals play an important role in this paper. We now recall their definition and basic properties.

\bg{Def}\cite{APT}\label{D-strong-idempotent}
A two-sided ideal $I$ of $\Lambda$ is called a strong idempotent ideal
if restriction of scalars induces isomorphisms
\[
\mr{Ext}^{i}_{\Lambda/I}(X,Y)\longrightarrow
\mr{Ext}^{i}_{\Lambda}(X,Y)
\qquad(i\geq0)
\]
for all $X,Y\in \Lambda/I\text{-}\mr{mod}$.
\ed{Def}

The following proposition gives two Ext-vanishing criteria for
strong idempotence.
\bg{Pro}\citep[Proposition~1.3]{APT}\label{Pro2.5}
Let $I$ be a two-sided ideal of an Artin algebra $\Lambda$. The following statements are equivalent.
\begin{enumerate}
\item[(a)] $I$ is a strong idempotent ideal.
\item[(b)] $\mr{Ext}^{i}_{\Lambda}(\Lambda/I,Y)=0$ for every
$Y\in\Lambda/I\text{-}\mr{mod}$ and every $i\geq1$.
\item[(c)] $\mr{Ext}^{i}_{\Lambda}(\Lambda/I,E)=0$ for every
injective object $E$ of $\Lambda/I\text{-}\mr{mod}$ and every $i\geq1$.
\end{enumerate}
\ed{Pro}

To study the relationship between \(\Lambda\) and its corner algebra $e\Lambda e$, Auslander, Platzeck, and Todorov \cite{APT} introduced the classes \(\mathbb{P}_k\) and \(\mathbb{I}_k\). We now recall their definitions. Let $e$ be an idempotent element of
$\Lambda$. For $k\in\mathbb N$, let $\mathbb P_k$ be the full
subcategory of $\Lambda\text{-}\mr{mod}$ consisting of modules
admitting a projective resolution $\cdots\to P_1\to P_0\to X\to0$
with $P_i\in\mr{add}(\Lambda e)$ for $0\leq i\leq k$.
The subcategory $\mathbb P_\infty$ consists of modules admitting
such a resolution with every term in $\mr{add}(\Lambda e)$.
Dually, $\mathbb I_k$ is the full subcategory of
$\Lambda\text{-}\mr{mod}$ consisting of modules admitting an
injective coresolution whose first $k+1$ terms belong to
$\mr{add}D(e\Lambda)$; for $\mathbb I_\infty$, every term is
required to belong to this subcategory.

To relate the homological properties of $\Lambda$ to those of
its quotient algebra and corner algebra, we recall the
\textbf{gluing sequences} established in
\citep[Proposition~3.3]{GLP}. Let $I=\Lambda e\Lambda$ be an idempotent ideal of
$\Lambda$, and set $\Gamma=\Lambda/I$.
For left $\Lambda$-modules $M$ and $X$,
we denote by $\tau_M(X)$ the trace of $M$ in $X$.
In particular,
\[
\tau_I(X)=IX=\Lambda eX,
\qquad
\tau_\Gamma(X)=\{x\in X\mid Ix=0\}.
\]
Thus, both $X/\tau_I(X)$ and $\tau_\Gamma(X)$
are naturally left $\Gamma$-modules. Following \citep[Section~5]{APT} and \cite{GLP},
we consider the perpendicular subcategories
\[
\mathbb T=\{X\in\Lambda\text{-}\mr{mod}\mid
\mr{Ext}^{i}_{\Lambda}(\Gamma,X)=0
\text{ for all }i\geq1\},
\]
\[
\widetilde{\mathbb T}=\{X\in\Lambda\text{-}\mr{mod}\mid
\mr{Ext}^{i}_{\Lambda}(X,D\Gamma)=0
\text{ for all }i\geq1\}.
\]
We recall the following gluing properties of these subcategories.

\begin{Lem}\label{L4.2}
{\rm\citep[Proposition~3.3]{GLP}}
Let $I=\Lambda e\Lambda$ be a strongly idempotent
ideal of $\Lambda$. Then the following statements hold.
\begin{enumerate}
\item[(a)]
$\mathbb I_\infty=\mathbb I_0\cap\mathbb T$.
For every $X\in\mathbb T$, there is a
\textbf{gluing sequence}
\[
0\longrightarrow\tau_\Gamma(X)
\longrightarrow X
\longrightarrow X/\tau_\Gamma(X)
\longrightarrow0,
\]
where
$\tau_\Gamma(X)\in\Gamma\text{-}\mr{mod}$
and
$X/\tau_\Gamma(X)\in\mathbb I_\infty$.

\item[(b)]
$\mathbb P_\infty
=\mathbb P_0\cap\widetilde{\mathbb T}$.
For every $X\in\widetilde{\mathbb T}$,
there is a \textbf{gluing sequence}
\[
0\longrightarrow\tau_I(X)
\longrightarrow X
\longrightarrow X/\tau_I(X)
\longrightarrow0,
\]
where
$\tau_I(X)\in\mathbb P_\infty$
and
$X/\tau_I(X)\in\Gamma\text{-}\mr{mod}$.
\end{enumerate}
\end{Lem}

\subsection{Auslander bounds}

We recall the definitions of the Auslander bounds
used throughout this paper; see \cite{W,L}.

\begin{Def}
For $M,N\in\Lambda\text{-}\mr{mod}$, set
\[
p_\Lambda(M,N)
=
\sup\{n\geq0\mid
\mr{Ext}_\Lambda^n(M,N)\neq0\}.
\]
The Auslander bound of $M$ is defined by
\[
b_\Lambda(M)
=
\sup\bigl\{
p_\Lambda(M,N)\mid
N\in\Lambda\text{-}\mr{mod},\
p_\Lambda(M,N)<\infty
\bigr\}.
\]
The global Auslander bound and the finitistic Auslander
bound of $\Lambda$ are defined, respectively, by
\[
\mr{gb}(\Lambda)
=
\sup\{b_\Lambda(M)\mid
M\in\Lambda\text{-}\mr{mod}\}
\]
and
\[
\mr{fab}(\Lambda)
=
\sup\bigl\{
b_\Lambda(M)\mid
M\in\Lambda\text{-}\mr{mod},\
b_\Lambda(M)<\infty
\bigr\}.
\]
All suprema are taken in
$\mathbb{N}\cup\{\infty\}$, with the convention
$\sup\varnothing=0$.

We say that $\Lambda$ satisfies the Auslander Condition
if $b_\Lambda(M)<\infty$ for every
$M\in\Lambda\text{-}\mr{mod}$, and that it satisfies
the Uniform Auslander Condition if $\mr{gb}(\Lambda)<\infty$.
\end{Def}

If $\mr{pd}_\Lambda M<\infty$, then
$b_\Lambda(M)=\mr{pd}_\Lambda M$
\citep[Remark 2.2(3)]{W}. Consequently,
\[
\mr{findim}(\Lambda)
\leq\mr{fab}(\Lambda)
\leq\mr{gb}(\Lambda).
\]
Moreover, $\mr{fab}(\Lambda)=\mr{gb}(\Lambda)$ whenever
$\Lambda$ satisfies the Auslander Condition.

In particular, finiteness of the global Auslander bound implies finiteness
of the finitistic Auslander bound, but the converse fails.
Indeed, the finite-dimensional commutative self-injective
algebra $A$ constructed in \cite{JS} satisfies
$\mr{gb}(A)=\infty$, whereas $\mr{fab}(A)=0$ by
\citep[Proposition~3.1]{W}.

\section{Auslander bounds under recollements}

Throughout this section, we fix an idempotent \(e\in\Lambda\) and set \(I=\Lambda e\Lambda\). Our aim is to compare the global and finitistic Auslander bounds of \(\Lambda\) with those of the quotient algebra \(\Lambda/I\) and the corner algebra \(e\Lambda e\). We first establish comparisons of Auslander bounds for individual modules via restriction and induction functors.
\begin{Lem}\label{L4.1}
Let $e$ be an idempotent element of $\Lambda$ and \(I=\Lambda e\Lambda\).
\begin{enumerate}
\item[(a)] If $I$ is a strong idempotent ideal, then,
for every $M\in\Lambda/I\text{-}\mr{mod}$,
\[
b_{\Lambda/I}(M)\leq b_{\Lambda}(M).
\]

\item[(b)] If $\mathbb{P}_{1}=\mathbb{P}_{\infty}$, then, for every
$M\in e\Lambda e\text{-}\mr{mod}$,
\[
b_{e\Lambda e}(M)=b_{\Lambda}(\Lambda e\otimes_{e\Lambda e}M).
\]
\end{enumerate}
\end{Lem}

\begin{proof}
(a)   Let
$L\in\Lambda/I\text{-}\mr{mod}$ satisfy
$
p_{\Lambda/I}(M,L)<\infty.
$
Since $I$ is a strong idempotent ideal,
Definition \ref{D-strong-idempotent} gives
\[
\mr{Ext}_{\Lambda/I}^{i}(M,L)
 \simeq
\mr{Ext}_{\Lambda}^{i}(M,L)
\qquad\text{for all }i\geq0.
\]
Consequently,
$
p_{\Lambda/I}(M,L)
 =p_{\Lambda}(M,L)
 \leq b_{\Lambda}(M).
$
Taking the supremum over all such $L$, we obtain
$
b_{\Lambda/I}(M)\leq b_{\Lambda}(M).
$

(b)
Let $M\in e\Lambda e\text{-}\mr{mod}$ and set $\widetilde M
=\Lambda e\otimes_{e\Lambda e}M.$ 
Recall that the functor
\[
\mr{Hom}_\Lambda(\Lambda e,-):
\mathbb P_1
\longrightarrow e\Lambda e\text{-}\mr{mod}
\]
is an equivalence with quasi-inverse
$\Lambda e\otimes_{e\Lambda e}-$.
Thus $\widetilde M\in\mathbb P_1
=\mathbb P_\infty.$ By \citep[Theorem~3.2(c)]{APT}, for every
$Y\in\Lambda\text{-}\mr{mod}$ and every $i\geq0$,
there are natural isomorphisms
\[
\begin{aligned}
\mr{Ext}_\Lambda^i(\widetilde M,Y)
&\cong
\mr{Ext}_{e\Lambda e}^i
\bigl(
\mr{Hom}_\Lambda(\Lambda e,\widetilde M),
\mr{Hom}_\Lambda(\Lambda e,Y)
\bigr)\\
&\cong
\mr{Ext}_{e\Lambda e}^i(M,eY).
\end{aligned}
\]
Consequently, $p_\Lambda(\widetilde M,Y)=p_{e\Lambda e}(M,eY)$
whenever either side is finite. It follows that
\[
b_{\Lambda}(\widetilde{M})
 \leq b_{e\Lambda e}(M).
\]

Conversely, let $L\in e\Lambda e\text{-}\mr{mod}$ and put $\widetilde{L}=\Lambda e\otimes_{e\Lambda e}L.$
Since $\mr{Hom}_{\Lambda}(\Lambda e,\widetilde{L})
 \simeq L,$
the preceding Ext isomorphism gives
$
\mr{Ext}_{\Lambda}^{i}
 (\widetilde{M},\widetilde{L})
 \simeq
\mr{Ext}_{e\Lambda e}^{i}(M,L)
$
for every $i\geq0$. Hence, $p_{e\Lambda e}(M,L)
 =
p_{\Lambda}(\widetilde{M},\widetilde{L})$
whenever either side is finite. Therefore,
$
b_{e\Lambda e}(M)
 \leq b_{\Lambda}(\widetilde{M}).
$
Combining the two inequalities, we conclude that
$
b_{e\Lambda e}(M)
 =
b_{\Lambda}(\Lambda e\otimes_{e\Lambda e}M).
$
\end{proof}
We next refine the comparison in Lemma~\ref{L4.1}(a) by establishing an upper bound in the opposite direction.
\begin{Pro}\label{P-fab-quotient}
Suppose that $I=\Lambda e\Lambda$
is a strong idempotent ideal and $r=\mr{pd}_{\Lambda}(\Lambda/I)<\infty$.
Then, for every $N\in\Lambda/I\text{-}\mr{mod}$, regarded also as a
left $\Lambda$-module by restriction of scalars, one has
\[
b_{\Lambda/I}(N)\leq b_\Lambda(N)\leq r+b_{\Lambda/I}(N).
\]
In particular, $b_{\Lambda/I}(N)<\infty$ if and only if
$b_\Lambda(N)<\infty$, and
\[
\mr{fab}(\Lambda/I)\leq\mr{fab}(\Lambda).
\]
\end{Pro}
\begin{proof}
The first inequality follows from Lemma \ref{L4.1}(a).
For the second, it suffices to consider $t=b_{\Lambda/I}(N)<\infty$.
Let $L\in\Lambda\text{-}\mr{mod}$ satisfy $p_\Lambda(N,L)<\infty$,
and put $Z=\Omega_\Lambda^{-r}L$, with $Z=L$ when $r=0$.
For every $j\geq1$, dimension shifting gives
\[
\mr{Ext}_\Lambda^j(\Lambda/I,Z)
\simeq\mr{Ext}_\Lambda^{j+r}(\Lambda/I,L)=0.
\]
Thus $Z\in\mathbb T$, and Lemma \ref{L4.2}(a) yields a gluing sequence
\[
0\longrightarrow U=\tau_{\Lambda/I}(Z)\longrightarrow Z
\longrightarrow V\longrightarrow0,
\]
where $U\in\Lambda/I\text{-}\mr{mod}$ and $V\in\mathbb I_\infty$.
By \citep[Proposition~2.6]{APT},
$\mr{Ext}_\Lambda^j(N,V)=0$ for every $j\geq0$.
The long exact Ext sequence, strong idempotence and dimension shifting
therefore give
\[
\mr{Ext}_{\Lambda/I}^i(N,U)
\simeq\mr{Ext}_\Lambda^i(N,Z)
\simeq\mr{Ext}_\Lambda^{i+r}(N,L)
\qquad (i\geq1).
\]
Consequently, $p_{\Lambda/I}(N,U)<\infty$, so $p_{\Lambda/I}(N,U)\leq t$.
For every $k>r+t$, taking $i=k-r>t$ in the preceding isomorphisms
gives $\mr{Ext}_\Lambda^k(N,L)=0$.
Hence $p_\Lambda(N,L)\leq r+t$.
Taking the supremum over all such $L$ proves
$b_\Lambda(N)\leq r+b_{\Lambda/I}(N)$.
In particular, if $b_{\Lambda/I}(N)<\infty$, then
\[
b_{\Lambda/I}(N)\leq b_\Lambda(N)\leq\mr{fab}(\Lambda).
\]
Taking the supremum over all such $N$ proves the final assertion.
\end{proof}

\begin{Lem}\label{L4.3}
Let $e$ be an idempotent element of $\Lambda$ and $I=\Lambda e\Lambda$.  Suppose that $\mr{pd}_{\Lambda}I=n<\infty.$
Then, for every $M\in\Lambda\text{-}\mr{mod}$, the following
statements hold.
\begin{enumerate}
\item[(a)] If
$\mr{gldim}(\Lambda/I)=l<\infty$ and $\mathbb{P}_{1}=\mathbb{P}_{\infty},$
then
\[
b_{\Lambda}(M)
 \leq
\max\{b_{e\Lambda e}(eM),n+l+2\}.
\]

\item[(b)] Suppose that
$\mr{pd}_{\Lambda^{op}}I=m<\infty.$
If $\mr{gldim}e\Lambda e=l<\infty$ and $I$ is a strong idempotent ideal, then
\[
b_{\Lambda}(M)
 \leq
m+n+2+
\max\left\{
l+1,
b_{\Lambda/I}
\left(
\frac{\Omega_{\Lambda}^{m+1}M}
{\tau_I(\Omega_{\Lambda}^{m+1}M)}
\right)
\right\}.
\]
\end{enumerate}
\end{Lem}
\begin{proof}
(a) Proposition \ref{P2.2}(e) gives an exact sequence
\[
0\to i_*(B)\to\Lambda e\otimes_{e\Lambda e}eM\to M
\to i_*((\Lambda/I)\otimes_\Lambda M)\to0,
\]
where $B\in\Lambda/I\text{-}\mr{mod}$.
Since $\mr{pd}_\Lambda(\Lambda/I)\leq n+1$
and $\mr{gldim}(\Lambda/I)=l$, both outer
nonzero terms have projective dimension at most $n+l+1$
over $\Lambda$. Splitting this sequence into two short exact
sequences and applying the long exact Ext sequences yields
 \[\mr{Ext}_\Lambda^k(M,-) \simeq \mr{Ext}_\Lambda^k(  \Lambda e\otimes_{e\Lambda e}eM,-) \qquad
\text{for all }k>  n+l+2.
 \]
 Put
$
r=b_{e\Lambda e}(eM).
$
If $r=\infty$, the asserted inequality is immediate. Suppose
that $r<\infty$. Since
$
\Lambda e\otimes_{e\Lambda e}eM
 \in\mathbb{P}_{1}
 =\mathbb{P}_{\infty},
$
it follows from \citep[Theorem~3.2(c)]{APT} that
$
\mr{Ext}_{\Lambda}^{i}
 \bigl(\Lambda e\otimes_{e\Lambda e}eM,L\bigr)
 \simeq
\mr{Ext}_{e\Lambda e}^{i}(eM,eL)
$
for every $L\in\Lambda\text{-}\mr{mod}$ and every $i\geq0$. Now let $L\in\Lambda\text{-}\mr{mod}$ satisfy
$p_{\Lambda}(M,L)<\infty$. Since
$
\mr{Ext}_{\Lambda}^{k}(M,L)
 \simeq
\mr{Ext}_{\Lambda}^{k}
 \bigl(\Lambda e\otimes_{e\Lambda e}eM,L\bigr)
$
for all $k>n+l+2$, it follows that
$p_{e\Lambda e}(eM,eL)<\infty$. Hence $p_{e\Lambda e}(eM,eL)\leq r.$
Therefore,
\[
\mr{Ext}_{\Lambda}^{k}(M,L)=0
\qquad
\text{for every }k>\max\{r,n+l+2\}.
\]
Thus
$
p_{\Lambda}(M,L)\leq\max\{r,n+l+2\}.
$
Since $L$ was arbitrary, we obtain
\[
b_{\Lambda}(M)
 \leq
\max\{b_{e\Lambda e}(eM),n+l+2\}.
\]

(b) Set
$ X=
\frac{\Omega_{\Lambda}^{m+1}M}{\tau_I(\Omega_{\Lambda}^{m+1}M)}$ and  $t=b_{\Lambda/I}(X)$. If $t=\infty$, the asserted inequality is immediate. Hence we
may assume that $t<\infty$. Since
$\mr{pd}_{\Lambda^{op}}(\Lambda/I)\leq m+1$,
dimension shifting gives
$\Omega_\Lambda^{m+1}M\in\widetilde{\mathbb T}$.
Lemma \ref{L4.2}(b) therefore yields a gluing sequence
\[
0 \to \tau_I( \Omega_{ \Lambda  }^{m+1} M )\to \Omega_{ \Lambda  }^{m+1} M \to X \to 0
\]
 with  $\tau_I ( \Omega_{ \Lambda  }^{m+1} M )\in \mathbb{P}_\infty $ and  $X\in \Lambda/I\text{-}\mr{mod}$.
Since $\mr{gldim}(e\Lambda e)=l$, we obtain
 \[
 \mr{pd}_{\Lambda}\tau_I ( \Omega_{ \Lambda  }^{m+1} M )=\mr{pd}_{e\Lambda e} \mr{Hom}_{\Lambda}(\Lambda e, \tau_I ( \Omega_{ \Lambda  }^{m+1} M ))\leq l
 \]
by \citep[Corollary~3.3(a)]{APT}. Hence
$
  \mr{Ext}_\Lambda^k(\Omega_{ \Lambda  }^{m+1} M ,-)\simeq \mr{Ext}_\Lambda^k( X,-) $ for all $k> l+1$.  Take any $L\in\Lambda\text{-}\mr{mod}$ with  $ p_{\Lambda}(M,L) < \infty$. Since $\mr{pd}_{\Lambda}I= n$, the $(n+1)$-th cosyzygy
$\Omega_{\Lambda}^{-n-1}L$ belongs to $\mathbb{T}$. By Lemma
\ref{L4.2}(a), there is a gluing sequence
\[
0\longrightarrow
\tau_{\Lambda/I}(\Omega_{\Lambda}^{-n-1}L)
\longrightarrow
\Omega_{\Lambda}^{-n-1}L
\longrightarrow
\frac{\Omega_{\Lambda}^{-n-1}L}
{\tau_{\Lambda/I}(\Omega_{\Lambda}^{-n-1}L)}
\longrightarrow 0,
\]
where
$
\tau_{\Lambda/I}(\Omega_{\Lambda}^{-n-1}L)
   \in\Lambda/I\text{-}\mr{mod}
$
and
$
\frac{\Omega_{\Lambda}^{-n-1}L}
{\tau_{\Lambda/I}(\Omega_{\Lambda}^{-n-1}L)}
   \in\mathbb{I}_{\infty}.$
 Set
\[
U=\tau_{\Lambda/I}(\Omega_{\Lambda}^{-n-1}L)
\qquad
\text{and}\qquad
V=
\frac{\Omega_{\Lambda}^{-n-1}L}
{\tau_{\Lambda/I}(\Omega_{\Lambda}^{-n-1}L)}.
\]
Thus, the preceding gluing sequence has the form
\[
0\longrightarrow U\longrightarrow
\Omega_{\Lambda}^{-n-1}L\longrightarrow V\longrightarrow0.
\]
Since $X\in\Lambda/I\text{-}\mr{mod}$ and
$V\in\mathbb{I}_{\infty}$, it follows from
\citep[Proposition~2.6]{APT} that
$
\mr{Ext}_{\Lambda}^{i}(X,V)=0$
 for all $i\geq0.$
Applying $\mr{Hom}_{\Lambda}(X,-)$ to the above
short exact sequence therefore yields
\[
\mr{Ext}_{\Lambda}^{i}(X,U)
 \simeq
\mr{Ext}_{\Lambda}^{i}
   (X,\Omega_{\Lambda}^{-n-1}L)
\qquad\text{for all }i\geq0.
\]
Since $I$ is a strong idempotent ideal and
$X,U\in\Lambda/I\text{-}\mr{mod}$,
Definition \ref{D-strong-idempotent} gives
\[
\mr{Ext}_{\Lambda}^{i}(X,U)
\simeq
\mr{Ext}_{\Lambda/I}^{i}(X,U)
\qquad(i\geq0).
\]
 We next show that $p_{\Lambda/I}(X,U)$ is finite.
Put
$
T=\tau_I(\Omega_{\Lambda}^{m+1}M).
$
As shown above, $\mr{pd}_{\Lambda}T\leq l$. From the gluing
sequence
\[
0\longrightarrow T\longrightarrow
\Omega_{\Lambda}^{m+1}M\longrightarrow X\longrightarrow0
\]
we obtain
$
\mr{Ext}_{\Lambda}^{i}
   (X,\Omega_{\Lambda}^{-n-1}L)
 \simeq
\mr{Ext}_{\Lambda}^{i}
   (\Omega_{\Lambda}^{m+1}M,\Omega_{\Lambda}^{-n-1}L)$
 for all $i>l+1.$
By dimension shifting, for every $i\geq1$,
\[
\mr{Ext}_{\Lambda}^{i}
   (\Omega_{\Lambda}^{m+1}M,\Omega_{\Lambda}^{-n-1}L)
 \simeq
\mr{Ext}_{\Lambda}^{i+m+n+2}(M,L).
\]
Since $p_{\Lambda}(M,L)<\infty$, the right-hand side vanishes
for all sufficiently large $i$. Consequently,
\[
p_{\Lambda/I}(X,U)<\infty.
\]
By the definition of $t$, we obtain
\[
p_{\Lambda/I}(X,U)\leq t.
\]
It follows that
\[
\mr{Ext}_{\Lambda}^{i}
   (\Omega_{\Lambda}^{m+1}M,\Omega_{\Lambda}^{-n-1}L)=0
\qquad\text{for all }i>\max\{l+1,t\}.
\]
Another application of dimension shifting gives
\[
\mr{Ext}_{\Lambda}^{q}(M,L)=0
\qquad\text{for all }
q>m+n+2+\max\{l+1,t\}.
\]
Therefore,
$
p_{\Lambda}(M,L)
 \leq m+n+2+\max\{l+1,t\}.
$
Since $L$ was arbitrary, we conclude that
\[
b_{\Lambda}(M)
 \leq m+n+2+\max\{l+1,t\}<\infty.
\]
\end{proof}
We now combine the preceding comparisons to obtain a reduction
theorem for the quotient algebra $\Lambda/I$.
Under suitable homological finiteness assumptions, the global
and finitistic Auslander bounds of $\Lambda$ can be controlled
by the corresponding bounds of the quotient. In particular,
these estimates yield an equivalence of finiteness for the
finitistic Auslander bounds.
\begin{Th}\label{T1}
Let $e$ be an idempotent element of $\Lambda$ and let $m,n,l\in\mathbb N$.
Suppose that $I$ is a strong idempotent ideal.
Then the following statements hold.
\begin{enumerate}
\item[(a)]  $\mr{gb}(\Lambda/I)\leq\mr{gb}(\Lambda).$
\item[(b)] If  $\mr{gldim}  e\Lambda e=l$, $\mr{pd}_{\Lambda}I= n$ and $\mr{pd}_{\Lambda^{op}}I=m$, then
\begin{align*}
\mr{gb}(\Lambda)&\leq m+n+2+
\max\left\{
l+1,
\mr{gb}(\Lambda/I)
\right\} \quad\text{and}\\
\mr{fab}(\Lambda)&\leq m+n+2+
\max\left\{
l+1,
\mr{fab}(\Lambda/I)
\right\}.
\end{align*}
Moreover, we have $\mr{fab}(\Lambda/I)
 \leq
\mr{fab}(\Lambda).$ Consequently,
\[
\mr{fab}(\Lambda)<\infty
\quad\Longleftrightarrow\quad
\mr{fab}(\Lambda/I)<\infty.
\]
\end{enumerate}
\end{Th}

\begin{proof}
(a) Let $M\in\Lambda/I\text{-}\mr{mod}$. By
Lemma \ref{L4.1}(a),
$b_{\Lambda/I}(M)\leq b_{\Lambda}(M)
 \leq\mr{gb}(\Lambda).$
Taking the supremum over all
$M\in\Lambda/I\text{-}\mr{mod}$, we obtain
$
\mr{gb}(\Lambda/I)
 \leq\mr{gb}(\Lambda).
$

(b) Suppose now that
$\mr{gldim}e\Lambda e=l<\infty,$ $\mr{pd}_{\Lambda}I=n<\infty,$
and
$\mr{pd}_{\Lambda^{op}}I=m<\infty.$
For every $M\in\Lambda\text{-}\mr{mod}$, put
\[
X_M=
\frac{\Omega_{\Lambda}^{m+1}M}
{\tau_I(\Omega_{\Lambda}^{m+1}M)}.
\]
By Lemma \ref{L4.3}(b),
$b_{\Lambda}(M)
 \leq
m+n+2+\max\{l+1,b_{\Lambda/I}(X_M)\}.$
Since
$b_{\Lambda/I}(X_M)
 \leq\mr{gb}(\Lambda/I),$
we have
$
b_{\Lambda}(M)
 \leq
m+n+2+
\max\left\{
l+1,\mr{gb}(\Lambda/I)
\right\}.
$
Taking the supremum over all
$M\in\Lambda\text{-}\mr{mod}$ gives
\[
\mr{gb}(\Lambda)
 \leq
m+n+2+
\max\left\{
l+1,\mr{gb}(\Lambda/I)
\right\}.
\]

We next prove the corresponding assertion for the finitistic
Auslander bound. Let $M\in\Lambda\text{-}\mr{mod}$ satisfy $b_{\Lambda}(M)<\infty,$
and write $X=\frac{\Omega_{\Lambda}^{m+1}M}
{\tau_I(\Omega_{\Lambda}^{m+1}M)}.$
We claim that
$
b_{\Lambda/I}(X )<\infty.
$
Indeed, let
$U\in\Lambda/I\text{-}\mr{mod}$ satisfy $p_{\Lambda/I}(X,U)<\infty.$
Since $I$ is a strong idempotent ideal,
Definition \ref{D-strong-idempotent} gives
$
p_{\Lambda}(X,U)
 =
p_{\Lambda/I}(X,U)<\infty.
$
Put $T=\tau_I(\Omega_{\Lambda}^{m+1}M).$
As in the proof of Lemma \ref{L4.3}(b), one has
$\mr{pd}_{\Lambda}T\leq l.$
Thus the exact sequence
\[
0\longrightarrow T\longrightarrow
\Omega_{\Lambda}^{m+1}M\longrightarrow X\longrightarrow0
\]
induces isomorphisms
$
\mr{Ext}_{\Lambda}^{i}(X,U)
 \simeq
\mr{Ext}_{\Lambda}^{i}
 (\Omega_{\Lambda}^{m+1}M,U)
 \simeq
\mr{Ext}_{\Lambda}^{i+m+1}(M,U)
$
for every $i>l+1$. It follows from
$p_{\Lambda}(X,U)<\infty$ that
$p_{\Lambda}(M,U)<\infty$. Hence $p_{\Lambda}(M,U)\leq b_{\Lambda}(M).$
The preceding Ext isomorphisms now imply
$
p_{\Lambda/I}(X,U)
 \leq
\max\{l+1,b_{\Lambda}(M)\}.
$
Taking the supremum over all such $U$, we obtain
$
b_{\Lambda/I}(X)
 \leq
\max\{l+1,b_{\Lambda}(M)\}<\infty.
$
Therefore,
$
b_{\Lambda/I}(X)
 \leq\mr{fab}(\Lambda/I).
$
Applying Lemma \ref{L4.3}(b) once more yields
$
b_{\Lambda}(M)
 \leq
m+n+2+
\max\left\{
l+1,\mr{fab}(\Lambda/I)
\right\}.
$
Taking the supremum over all $M$ satisfying
$b_{\Lambda}(M)<\infty$, we conclude that
\[
\mr{fab}(\Lambda)
 \leq
m+n+2+
\max\left\{
l+1,\mr{fab}(\Lambda/I)
\right\}.
\]

Finally, the exact sequence
$0\to I\to\Lambda\to\Lambda/I\to0$
gives $\mr{pd}_{\Lambda}(\Lambda/I)\leq n+1<\infty$.
Proposition \ref{P-fab-quotient} therefore yields
$
\mr{fab}(\Lambda/I)\leq\mr{fab}(\Lambda).
$
Together with the upper bound just proved, this establishes
the asserted equivalence of finiteness.
\end{proof}
We next turn to the corner algebra $e\Lambda e$ and obtain
the complementary reduction theorem.
\begin{Th}\label{T2}
Let $e$ be an idempotent element of $\Lambda$ and let $n,l\in\mathbb N$.
Suppose that $\mathbb P_1=\mathbb P_\infty$.
Then the following statements hold.
\begin{enumerate}
\item[(a)] $\mr{gb}(  e\Lambda e)\leq\mr{gb}(\Lambda)$  and $\mr{fab}(e\Lambda e)\leq\mr{fab}(\Lambda).$
\item[(b)]  If $\mr{gldim}(\Lambda/I)=l<\infty$, $\mr{pd}_{\Lambda}I= n$,  then
\begin{align*}
\mr{gb}(\Lambda)&\leq\max\{n+l+2,\mr{gb}(e\Lambda e)\}&\text{and} \qquad
\mr{fab}(\Lambda)&\leq\max\{n+l+2,\mr{fab}(e\Lambda e)\}.
\end{align*}
\end{enumerate}
\end{Th}

\begin{proof}
(a)  For every
$M\in e\Lambda e\text{-}\mr{mod}$, Lemma \ref{L4.1}(b) gives
$
b_{e\Lambda e}(M)=
b_{\Lambda}(\Lambda e\otimes_{e\Lambda e}M)
 \leq\mr{gb}(\Lambda).
$
Taking the supremum over all
$M\in e\Lambda e\text{-}\mr{mod}$ gives
$
\mr{gb}(e\Lambda e)\leq\mr{gb}(\Lambda).
$

We next prove the corresponding inequality for the finitistic
Auslander bound. If $b_{e\Lambda e}(M)<\infty$, then
$b_{\Lambda}(\Lambda e\otimes_{e\Lambda e}M)<\infty$
and
$
b_{e\Lambda e}(M)
 =
b_{\Lambda}(\Lambda e\otimes_{e\Lambda e}M)
 \leq\mr{fab}(\Lambda).
$
Taking the supremum over all $M$ satisfying
$b_{e\Lambda e}(M)<\infty$ gives
$
\mr{fab}(e\Lambda e)\leq\mr{fab}(\Lambda).
$

(b) Suppose further that
$\mr{gldim}(\Lambda/I)=l<\infty$ and $\mr{pd}_{\Lambda}I=n<\infty.$
By Lemma \ref{L4.3}(a), for every
$M\in\Lambda\text{-}\mr{mod}$,
$
b_{\Lambda}(M)
 \leq
\max\{n+l+2,b_{e\Lambda e}(eM)\}.
$
Since $b_{e\Lambda e}(eM)\leq\mr{gb}(e\Lambda e),$  taking the supremum over all $M$ gives
$
\mr{gb}(\Lambda)
 \leq
\max\{n+l+2,\mr{gb}(e\Lambda e)\}.
$

Finally, let $M\in\Lambda\text{-}\mr{mod}$ satisfy
$b_{\Lambda}(M)<\infty$. We first show that $b_{e\Lambda e}(eM)<\infty.$
Let $Y\in e\Lambda e\text{-}\mr{mod}$ satisfy $p_{e\Lambda e}(eM,Y)<\infty,$
and put $\widetilde{Y}=\Lambda e\otimes_{e\Lambda e}Y.$
Since $\mathbb{P}_{1}=\mathbb{P}_{\infty}$,
\[
\mr{Ext}_{\Lambda}^{i}
 (\Lambda e\otimes_{e\Lambda e}eM,\widetilde{Y})
 \simeq
\mr{Ext}_{e\Lambda e}^{i}(eM,Y)
\]
for every $i\geq0$. Moreover, the proof of Lemma
\ref{L4.3}(a) gives
\[
\mr{Ext}_{\Lambda}^{i}(M,\widetilde{Y})
 \simeq
\mr{Ext}_{\Lambda}^{i}
 (\Lambda e\otimes_{e\Lambda e}eM,\widetilde{Y})
\]
for every $i>n+l+2$. It follows from $p_{e\Lambda e}(eM,Y)<\infty$ that
$p_{\Lambda}(M,\widetilde{Y})<\infty$. Therefore, $
p_{\Lambda}(M,\widetilde{Y})\leq b_{\Lambda}(M),$
and the preceding isomorphisms imply
$
p_{e\Lambda e}(eM,Y)
 \leq
\max\{n+l+2,b_{\Lambda}(M)\}.
$
Taking the supremum over all such $Y$, we obtain
$
b_{e\Lambda e}(eM)
 \leq
\max\{n+l+2,b_{\Lambda}(M)\}<\infty.
$
Hence $b_{e\Lambda e}(eM)\leq\mr{fab}(e\Lambda e).$
Applying Lemma \ref{L4.3}(a), we get
$
b_{\Lambda}(M)
 \leq
\max\{n+l+2,\mr{fab}(e\Lambda e)\}.
$
Taking the supremum over all $M$ for which
$b_{\Lambda}(M)<\infty$ gives
\[
\mr{fab}(\Lambda)
 \leq
\max\{n+l+2,\mr{fab}(e\Lambda e)\}.
\]
\end{proof}
The following examples illustrate the two complementary reduction
mechanisms in Theorems \ref{T1} and \ref{T2}.
\begin{Exm}\label{E-T1-nontriangular}
Let
$\Lambda=\Bbbk Q/\langle ab,a\varepsilon,\varepsilon b,\varepsilon^2\rangle,$
where \(Q\) is the quiver \[
\qquad
\xymatrix{
1 \ar@/^/[r]^a
&
2 \ar@/^/[l]^b
  \ar@(ur,dr)[]^{\varepsilon}.
}
\]
Throughout this example, paths are composed
from left to right; thus $ab$ means first traversing $a$
and then $b$. Put
$e=e_1$ and $I=\Lambda e\Lambda$.
Then $e\Lambda e\simeq\Bbbk,$ and $
\Lambda/I
\simeq D:=\Bbbk[\varepsilon]/(\varepsilon^2).
$
Since $I
=\langle e_1,a,b,ba\rangle_{\Bbbk},$ we have $
{}_{\Lambda}I
\simeq(\Lambda e)^{\oplus2},$ and $I
_{\Lambda}\simeq(e\Lambda)^{\oplus2}.
$
Thus $I
$ is projective on both sides. By \citep[Proposition~5.2]{APT},   $I $ is a  strong idempotent ideal and  $\mathbb{P}_{1}=\mathbb{P}_{\infty}$.
We claim that $\mr{gb}(D)=\mr{fab}(D)=0.$
In fact, let $S=D/(\varepsilon)$. Every finitely generated
$D$-module is a direct sum of copies of $D$ and $S$.
The periodic projective resolution
\[
\cdots\xrightarrow{\cdot\varepsilon}D
\xrightarrow{\cdot\varepsilon}D
\xrightarrow{\cdot\varepsilon}D
\longrightarrow S\longrightarrow0
\]
gives
$
\mr{Ext}_D^i(S,S)\simeq\Bbbk$ and
$\mr{Ext}_D^i(S,D)=0
 $ $(i\geq1)$.
Together with projectivity of $D$, these identities
show that eventual vanishing of
$\mr{Ext}_D^i(M,N)$ implies vanishing for
every $i\geq1$. Consequently,
\[
\mr{gb}(D)=\mr{fab}(D)=0.
\]
We can now apply Theorem \ref{T1}(b) with
$m=n=l=0$. It gives
\[
\mr{gb}(\Lambda)
\leq 2+\max\{1,\mr{gb}(D)\}=3
\]
and
\[
\mr{fab}(\Lambda)
\leq 2+\max\{1,\mr{fab}(D)\}=3.
\]
Clearly, $\mr{gldim}\Lambda=\infty$.
 The example therefore illustrates reduction of
Auslander bounds through the quotient algebra,
without finite global dimension of either
$\Lambda/I $ or $\Lambda$.
\end{Exm}

\begin{Exm}\label{E-T2-nontriangular}
Let
$\Lambda=\Bbbk Q/\langle x^2,ab\rangle,$
where \(Q\) is the quiver
\[
\xymatrix{
1 \ar@/^/[r]^a
  \ar@(ul,dl)[]_{x}
&
2 \ar@/^/[l]^b.
}
\]
 Again,
paths are composed from left to right. Put
$e=e_1$ and $I=\Lambda e\Lambda$.
Then
$
e\Lambda e\simeq D:=\Bbbk[x]/(x^2)$  and $
\Lambda/I \simeq\Bbbk.$
Since $I
=\langle e_1,x,a,xa,b,bx,ba,bxa\rangle_{\Bbbk},$ we have $
{}_{\Lambda}I
\simeq(\Lambda e)^{\oplus2},$ and $I
_{\Lambda}\simeq(e\Lambda)^{\oplus2}.
$
Thus $I
$ is projective on both sides. As in Example \ref{E-T1-nontriangular},
we have $\mathbb P_1=\mathbb P_\infty$ and
$\mr{gb}(D)=\mr{fab}(D)=0$.
Applying Theorem  \ref{T2}(b)  with
$n=l=0$, we have
\[
\mr{gb}(\Lambda)
\leq \max\{2,\mr{gb}(D)\}=2 \qquad \text{and}\qquad \mr{fab}(\Lambda)
\leq \max\{2,\mr{fab}(D)\}=2.
\]
We now show that both upper bounds are attained. Put $P_i=\Lambda e_i$ and
$S_i=P_i/\mr{rad}P_i$.
Right multiplication by $a$ gives an exact sequence
\[
0\longrightarrow P_1
\xrightarrow{\cdot a}P_2
\longrightarrow S_2\longrightarrow0.
\]
Since $\mr{Hom}_{\Lambda}(P_2,S_1)=0$ and
$\mr{Hom}_{\Lambda}(P_1,S_1)\simeq\Bbbk,$
we obtain
\[
\mr{Ext}_{\Lambda}^1(S_2,S_1)
\simeq\Bbbk.
\]
In particular, $\mr{pd}_{\Lambda}S_2=1$.  Set $X=P_1/\mr{soc}P_1$. The exact sequence
\[
0\longrightarrow S_2^{\oplus2}
\longrightarrow P_1\longrightarrow X
\longrightarrow0
\]
shows that $\mr{pd}_{\Lambda}X=2$, while
dimension shifting gives
\[
\mr{Ext}_{\Lambda}^2(X,S_1)
\simeq
\mr{Ext}_{\Lambda}^1
(S_2^{\oplus2},S_1)
\simeq\Bbbk^2.
\]
Therefore
$\mr{pd}_{\Lambda}X=b_{\Lambda}(X)=2.$
Combining this with the upper bounds proves
\[
\mr{gb}(\Lambda)=\mr{fab}(\Lambda)=2.
\]
\end{Exm}

Applying Theorems \ref{T1} and \ref{T2} to upper triangular matrix algebras,
we obtain the following corollary.

\begin{Cor}\label{C-upper-triangular-A}
Let $\Lambda=
\begin{pmatrix}
A&{}_A U_B\\
0&B
\end{pmatrix}.$
Then
\begin{enumerate}
\item[(a)]$
\max\{\mr{gb}( A),\mr{gb}( B)\}\leq\mr{gb}(\Lambda) $  and  $\mr{fab}(A)\leq\mr{fab}(\Lambda).
$
Moreover, if $\mr{pd}_{B^{op}}U= 0$ or $\mr{gldim}A=l<\infty$, then
\[
\max\{\mr{fab}( A),\mr{fab}(B)\}\leq\mr{fab}(\Lambda).
\]
\item[(b)]  If $\mr{gldim}A=l<\infty$,  then
\begin{align*}
\mr{gb}(\Lambda)&\leq\max\{2l+3,\mr{gb}(B)+l+2\}&\text{and} \qquad
\mr{fab}(\Lambda)&\leq\max\{2l+3,\mr{fab}(B)+l+2\}.
\end{align*}
\item[(c)] If $\mr{gldim}B=l<\infty$,  then
\begin{align*}
\mr{gb}(\Lambda)&\leq\max\{2l+3,\mr{gb}(A)+l+2\}&\text{and} \qquad
\mr{fab}(\Lambda)&\leq\max\{2l+3,\mr{fab}(A)+l+2\}.
\end{align*}
\end{enumerate}
\end{Cor}

\begin{proof}
Let
$e_1=\begin{pmatrix}1&0\\0&0\end{pmatrix}$ and $e_2=\begin{pmatrix}0&0\\0&1\end{pmatrix}.$ Take $e=e_1$, and put $I=\Lambda e\Lambda$. Then
$ e\Lambda e \simeq A,$ $ \Lambda/I\simeq B $ and  $ I\simeq e_1\Lambda.
$
The ideal $I=e_1\Lambda$ is projective as a right
$\Lambda$-module. By \citep[Proposition~5.2]{APT}, it is a strong
idempotent ideal and $\mathbb P_1=\mathbb P_\infty$ for $e=e_1$.
Theorems \ref{T1}(a) and \ref{T2}(a) therefore yield
 \[
\max\{\mr{gb}( A),\mr{gb}( B)\}\leq\mr{gb}(\Lambda) \qquad  \text{and} \qquad \mr{fab}(A)\leq\mr{fab}(\Lambda).
\]
This proves the first part of (a). Suppose that $\mr{gldim}A=l<\infty$. 
We have
$
\mr{pd}_{\Lambda}(I)=\mr{pd}_{A}U\leq l.
$
 Since $\mr{pd}_{\Lambda^{op}}I=0$, Theorem \ref{T1}(b) yields
\begin{align*}
\mr{gb}(\Lambda)&\leq\max\{2l+3,\mr{gb}(B)+l+2\},\\
\mr{fab}(\Lambda)&\leq\max\{2l+3,\mr{fab}(B)+l+2\},
\end{align*}
and $\mr{fab}(B)\leq\mr{fab}(\Lambda)$.
This proves (b) and the second assertion of (a) when $A$ has
finite global dimension.

Now take $e=e_2$, with $I=\Lambda e\Lambda$. Then
$e\Lambda e\simeq B,$
$\Lambda/I\simeq A$ and
$I=\Lambda e.$
As a right $\Lambda$-module, $I$ is the module
supported on the $B$-component with underlying $B$-module $U\oplus B$.
Thus $\mr{pd}_{\Lambda^{op}}(I)=\mr{pd}_{B^{op}}(U\oplus B)$.
If $\mr{pd}_{B^{op}}U=0$, this ideal is projective on the right.
By \citep[Proposition~5.2]{APT}, $\mathbb P_1=\mathbb P_\infty$
for $e=e_2$, and Theorem \ref{T2}(a) yields
$
\mr{fab}(B)\leq\mr{fab}(\Lambda).
$
This proves the remaining assertion in (a).
Finally, we prove (c). Since $I=\Lambda e$ is
projective as a left $\Lambda$-module, we have
$
\mr{pd}_{\Lambda}(I)=0.
$
The exact sequence of left $\Lambda$-modules \[ 0\longrightarrow\Lambda e\longrightarrow\Lambda \longrightarrow\Lambda/I\longrightarrow0 \] is a projective resolution of $\Lambda/I$. For every $Y\in\Lambda/I\text{-}\mr{mod}$,
$\mr{Hom}_{\Lambda}(\Lambda e,Y)\simeq eY=0. $ Hence \[\mr{Ext}_{\Lambda}^{i} (\Lambda/I,Y)=0 \qquad\text{for all }i\geq1, \] so Proposition \ref{Pro2.5}(b) shows that $I$ is a strong idempotent ideal.
Now, suppose that $\mr{gldim}B=l<\infty$.
Since
$
\mr{pd}_{\Lambda^{op}}(I)
=\mr{pd}_{B^{op}}(U\oplus B)\leq l,
$
(c) follows  from Theorem \ref{T1}(b), with
$n=0$ and $m\leq l$.
\end{proof}


%
The preceding corollary follows from the general recollement estimates.
The triangular structure and Proposition \ref{P-fab-quotient} yield
the following sharper bounds.

\begin{Pro}\label{P-upper-triangular-sharp}
Let $\Lambda=\begin{pmatrix}A&{}_A U_B\\0&B\end{pmatrix}$.
Then the following statements hold.
\begin{enumerate}
\item[(a)] Suppose that $\mr{gldim}A=l<\infty$, and put
$d=\mr{pd}_A U$, with $d=0$ if $U=0$. Then
\begin{align*}
\mr{gb}(\Lambda)
&\leq\max\{l+1,d+1+\mr{gb}(B)\}
\leq l+1+\mr{gb}(B),\\
\mr{fab}(\Lambda)
&\leq\max\{l+1,d+1+\mr{fab}(B)\}
\leq l+1+\mr{fab}(B).
\end{align*}
\item[(b)] If $\mr{gldim}B=l<\infty$, then
\[
\mr{gb}(\Lambda)\leq l+1+\mr{gb}(A)
\qquad\text{and}\qquad
\mr{fab}(\Lambda)\leq l+1+\mr{fab}(A).
\]
\end{enumerate}
Consequently, in (a), finiteness of either bound for $\Lambda$
is equivalent to finiteness of the corresponding bound for $B$;
in (b), it is equivalent to finiteness of the corresponding bound for $A$.
\end{Pro}
\begin{proof}
We write a left $\Lambda$-module as
$
M=\binom{X}{Y}_{\phi},
$
where $X\in A\text{-}\mr{mod}$,
$Y\in B\text{-}\mr{mod}$, and
$\phi:U\otimes_B Y\to X$ is $A$-linear.
We omit the subscript when the structure map is clear.
The functor $F:A\text{-}\mr{mod}\longrightarrow
\Lambda\text{-}\mr{mod},$
 $F(X)=\binom{X}{0},$
is exact and preserves projectives. Applying $F$ to a projective resolution gives
\[
\mr{Ext}_{\Lambda}^i(\binom{X}{0},\binom{X'}{Y'})
\simeq\mr{Ext}_A^i(X,X')\qquad(i\geq0).
\]
Every $A$-module occurs as the first component of $\binom{X'}{0}$, so
$
b_\Lambda(\binom{X}{0})=b_A(X).$ Likewise, the exact functor $e_2(-)$ preserves projectives, since
$e_2\Lambda e_1=0$ and $e_2\Lambda e_2=B$.
Applying it to a projective resolution, and using the natural
Hom identification, yields
\[
\mr{Ext}_\Lambda^i(M,\binom{0}{V})
\simeq\mr{Ext}_B^i(Y,V)\qquad(i\geq0).
\]
In particular, $b_B(Y)\leq b_\Lambda(M).$

(a) Take $e=e_1$, and put $I=\Lambda e\Lambda$. The ideal $I$ is strong idempotent, as in
the proof of Corollary \ref{C-upper-triangular-A}.
As a left $\Lambda$-module, $I=\binom{A\oplus U}{0}$.
Hence $\mr{pd}_\Lambda I\leq d$, and
$\mr{pd}_\Lambda B\leq d+1$, where $\Lambda/I\simeq B$.
Proposition \ref{P-fab-quotient} gives
\[
b_\Lambda(\binom{0}{Y})\leq d+1+b_B(Y).
\]
The canonical exact sequence
\[
0\longrightarrow\binom{X}{0}
\longrightarrow\binom{X}{Y}_{\phi}
\longrightarrow\binom{0}{Y}
\longrightarrow0
\]
and $\mr{pd}_\Lambda\binom{X}{0}=\mr{pd}_A X\leq l$ imply
\[
\mr{Ext}_\Lambda^i(\binom{0}{Y},L)
\simeq\mr{Ext}_\Lambda^i(\binom{X}{Y}_{\phi},L)\qquad(i>l+1)
\]
for every left $\Lambda$-module $L$.
Thus, whenever $p_\Lambda(M,L)<\infty$, one also has
$p_\Lambda(\binom{0}{Y},L)<\infty$. It follows that
\[
b_\Lambda(M)
\leq\max\{l+1,b_\Lambda(\binom{0}{Y})\}
\leq\max\{l+1,d+1+b_B(Y)\}.
\]
Taking the supremum over all $M$ proves the global bound.
If $b_\Lambda(M)<\infty$, the inequality
$b_B(Y)\leq b_\Lambda(M)$ proved above guarantees
$b_B(Y)<\infty$; hence $b_B(Y)\leq\mr{fab}(B)$.
Taking the supremum over these $M$ proves the finitistic bound.
The remaining inequalities follow from $d\leq l$.

(b) Put $s=l+1$. We first show that a projective resolution of $M$
can be chosen so that its $s$th syzygy is $\binom{X_s}{0}$
for some $A$-module $X_s$.
Indeed, choose projective epimorphisms $Q\to Y$ over $B$ and
$P\to X$ over $A$. They determine a projective epimorphism
\[
\binom{P}{0}\oplus \binom{U\otimes_B Q}{Q}_{\mathrm{id}}\twoheadrightarrow
\binom{X}{Y}_{\phi}.
\]
 Repeating this construction along a projective
$B$-resolution of $Y$, its $l$th syzygy has projective second
component. At the next step, choose the map on that component
to be the identity. Its kernel then has zero second component,
as asserted.

Write $K=\Omega_\Lambda^sM=\binom{X_s}{0}$.
Dimension shifting gives
$\mr{Ext}_\Lambda^i(K,L)
\simeq\mr{Ext}_\Lambda^{i+s}(M,L),$ $i\geq1.$
In particular, eventual vanishing for either pair is equivalent
to eventual vanishing for the other, and
\[
b_\Lambda(K)\leq b_\Lambda(M),
\qquad
b_\Lambda(M)\leq s+b_\Lambda(K).
\]
Since $b_\Lambda(K)=b_A(X_s)$, we obtain
$b_\Lambda(M)\leq l+1+b_A(X_s)$.
Taking the supremum over all $M$ proves the global bound.
If $b_\Lambda(M)<\infty$, then
$b_A(X_s)=b_\Lambda(K)<\infty$, so the same argument with
$\mr{fab}(A)$ proves the finitistic bound.

Finally, the corresponding lower bounds in
Corollary \ref{C-upper-triangular-A}(a) give the asserted
equivalences of finiteness.
\end{proof}

We conclude with an example illustrating sharpness and the applicability
of the estimates beyond finite global dimension.

\begin{Exm}\label{E-triangular-infinite-gldim}
Put $D=\Bbbk[\varepsilon]/(\varepsilon^2)$, which has infinite
global dimension. By Example \ref{E-T1-nontriangular},
\[
\mr{gb}(D)=\mr{fab}(D)=0.
\]
Consider the two upper triangular matrix algebras
\[
\Lambda_1=
\begin{pmatrix}\Bbbk&{}_{\Bbbk}D_D\\0&D\end{pmatrix},
\qquad
\Lambda_2=
\begin{pmatrix}D&{}_D D_{\Bbbk}\\0&\Bbbk\end{pmatrix}.
\]
Proposition \ref{P-upper-triangular-sharp}(a) and (b), respectively,
give $\mr{gb}(\Lambda_j)\leq1$ and
$\mr{fab}(\Lambda_j)\leq1$ for $j=1,2$.
To see that these estimates are exact, write
$\Lambda_j=\begin{pmatrix}A&U\\0&B\end{pmatrix}$ in each case.
The left $A$-module $U$ is nonzero and projective, and there is
an exact sequence
\[
0\longrightarrow \binom{U}{0}\longrightarrow
\Lambda_j e_2=\binom{U}{B}_{ \mathrm{id}}
\longrightarrow\binom{0}{B}\longrightarrow0.
\]
The first two terms are projective, but this sequence does not split:
a section with second component $\mathrm{id}_B$ would violate
the compatibility with the nonzero structure map
$U\otimes_B B\to U$.
It therefore represents a nonzero element of
$\mr{Ext}_{\Lambda_j}^1(\binom{0}{B},\binom{U}{0})$.
Thus $b_{\Lambda_j}(\binom{0}{B})=1$, and
\[
\mr{gb}(\Lambda_j)=\mr{fab}(\Lambda_j)=1
\qquad(j=1,2).
\]
Nevertheless, both algebras have infinite global dimension. These examples show that the new estimates apply even when
the other corner and the triangular matrix algebra have infinite
global dimension. They also improve the preceding corollary's
upper bound $3$ to the exact value $1$.
\end{Exm}

\section{Auslander bounds and arrow removals}
In this section, we investigate whether finiteness of the global
and finitistic Auslander bounds is preserved under arrow removal.
We first study these questions in the more general framework of
cleft extensions and then apply the resulting comparison
theorems to arrow removal. We begin by recalling the terminology
of cleft extensions used below. We denote the faithful exact functor
by $\mathsf{e}$ to distinguish it from the idempotent $e$ used
in the preceding sections.
\begin{Def}\label{def:cleft-extension}
{\rm\cite[Definition~2.1]{Beligiannis2000}}
A \emph{cleft extension} of an abelian category $\cB$ consists of a
diagram of functors
\[
\xymatrix@C=3.7pc{
\cB\ar[r]^-i&
\cA\ar[r]^-{\mathsf{e}}&
\cB\ar@<-1.8ex>[l]_-l
}
\]
such that $\mathsf{e}$ is faithful and exact, $(l,\mathsf{e})$ is an adjoint pair, and
there is a natural isomorphism
$
\mathsf{e} i \cong \id_{\cB}.
$
We denote such a cleft extension by $(\cB,\cA,i,\mathsf{e},l)$.
\end{Def}

 Let
$({\mathcal B},{\mathcal A},i,\mathsf{e},l)$
be a cleft extension of categories of finitely generated modules.
For $t\in\mathbb N$, the exact functor $\mathsf{e}$ is called a
$t$-eventually homological isomorphism if the natural maps
induced by $\mathsf{e}$ give isomorphisms
\begin{equation}\label{cleft-ext-isomorphism}
\mr{Ext}^{q}_{\mathcal A}(X,Y)\simeq
\mr{Ext}^{q}_{\mathcal B}(\mathsf{e}X,\mathsf{e}Y)
\end{equation}
for all $X,Y\in\mathcal A$ and all $q>t$.
The invariants $b$, $\mr{gb}$ and $\mr{fab}$ for these module
categories are defined as in Section 2.
\begin{Pro}\label{P-cleft-bound}
Let $({\mathcal B},{\mathcal A},i,\mathsf{e},l)$ be a cleft extension of
categories of finitely generated modules. Suppose that $\mathsf{e}$ is a
$t$-eventually homological isomorphism for some $t\in\mathbb N$.
Then the following statements hold.
\begin{enumerate}
\item[(a)]
For every $X\in{\mathcal A}$,
$$b_{\mathcal A}(X)\leq\max\{t,b_{\mathcal B}(\mathsf{e}X)\}
\quad\text{and}\quad
b_{\mathcal B}(\mathsf{e}X)\leq\max\{t,b_{\mathcal A}(X)\}.$$
Consequently, $b_{\mathcal A}(X)<\infty$ if and only if $b_{\mathcal B}(\mathsf{e}X)<\infty$.
\item[(b)]
\begin{align*}
\mr{gb}({\mathcal A})&\leq\max\{t,\mr{gb}({\mathcal B})\},&
\mr{gb}({\mathcal B})&\leq\max\{t,\mr{gb}({\mathcal A})\},\\
\mr{fab}({\mathcal A})&\leq\max\{t,\mr{fab}({\mathcal B})\},&
\mr{fab}({\mathcal B})&\leq\max\{t,\mr{fab}({\mathcal A})\}.
\end{align*}
Consequently, $\mr{gb}({\mathcal A})<\infty$ if and only if $\mr{gb}({\mathcal B})<\infty$, and $\mr{fab}({\mathcal A})<\infty$ if and only if $\mr{fab}({\mathcal B})<\infty$.
\end{enumerate}
\end{Pro}
\begin{proof}
(a) Let $Y\in{\mathcal A}$ and suppose that $p_{\mathcal A}(X,Y)<\infty$. By \eqref{cleft-ext-isomorphism}, $p_{\mathcal B}(\mathsf{e}X,\mathsf{e}Y)<\infty$, whence
$$p_{\mathcal A}(X,Y)\leq\max\{t,b_{\mathcal B}(\mathsf{e}X)\}.$$
Taking the supremum over all such $Y$ proves the first inequality. Conversely, let $V\in{\mathcal B}$ and suppose that $p_{\mathcal B}(\mathsf{e}X,V)<\infty$. Since $\mathsf{e}i\simeq\operatorname{id}_{\mathcal B}$, the functor $\mathsf{e}$ is essentially surjective.  There is a $Y\in{\mathcal A}$ with $\mathsf{e}Y\simeq V$. Equation \eqref{cleft-ext-isomorphism} implies $p_{\mathcal A}(X,Y)<\infty$, and hence
$$p_{\mathcal B}(\mathsf{e}X,V)\leq\max\{t,b_{\mathcal A}(X)\}.$$
Taking the supremum over all such $V$ gives the second inequality.

(b) The assertions for the global bounds follow by taking suprema in (a) and using the essential surjectivity of $\mathsf{e}$. (a) also identifies the objects having finite individual Auslander bound; taking suprema over precisely these objects proves the assertions for the finitistic bounds.

\end{proof}

The arrow-removal operation and its homological properties have been
studied in \cite{EPS,GPS2021}. We work with left modules,
whereas these papers use right modules.
Let $Q$ be a finite quiver and let $\Lambda=\Bbbk Q/I$ be an admissible quotient of the path algebra $\Bbbk Q$. Suppose that $Q$ contains arrows
$$a_i: v_{e_i}\longrightarrow v_{f_i},\qquad 1\leq i\leq s,$$
which do not occur in a set of minimal generators of $I$ and satisfy
$$\mr{Hom}_{\Lambda}(\Lambda f
_i,\Lambda e_j)=0
  \qquad\text{for all }1\leq i,j\leq s.$$
Set
$\Gamma=\Lambda/\Lambda\{\overline a_i\}_{i=1}^s\Lambda,$
where $\overline a_i=a_i+I$. The natural projection $\Lambda\to\Gamma$ is called an arrow removal.  Write $\pi:\Lambda\to\Gamma$ for the canonical projection.
The arrow-removal construction admits an algebra section
$\sigma:\Gamma\to\Lambda$ satisfying
$\pi\sigma=\operatorname{id}_{\Gamma}$.
In the following diagram, the $\Lambda$-module structures on
$\Gamma$ are induced by $\pi$, whereas the $\Gamma$-module
structures on $\Lambda$ are induced by $\sigma$.
Thus $i$ and $\mathsf{e}$ are naturally isomorphic to restriction of
scalars along $\pi$ and $\sigma$, respectively.
 The arrow-removal construction gives the following diagram of functors:
\[
\begin{tikzcd}[column sep=5.5em]
\Gamma\text{-}\mr{mod}
\arrow[bend left=0]{r}{i=\mr{Hom}_\Gamma({}_\Gamma \Gamma_\Lambda,-)}
&
\Lambda\text{-}\mr{mod}
\arrow[bend left=0]{r}{\mathsf{e}=\mr{Hom}_\Lambda({}_\Lambda \Lambda_\Gamma,-)}
\arrow[bend left=40]{l}{p=\mr{Hom}_\Lambda({}_ \Lambda\Gamma_\Gamma,-)}
\arrow[bend right=40]{l}[above]{q={}_\Gamma \Gamma_\Lambda\otimes_\Lambda- }
&
\Gamma\text{-}\mr{mod}.
\arrow[bend left=40]{l}{r=\mr{Hom}_\Gamma({}_\Gamma \Lambda_\Lambda,-)}
\arrow[bend right=40]{l}[above]{l={}_\Lambda\Lambda_\Gamma \otimes_\Gamma- }
\end{tikzcd}
\]
We now show that the associated module categories form a cleft
extension satisfying the hypotheses of Proposition \ref{P-cleft-bound}.
To apply the right-module results of \cite{GPS2021,EPS},
we pass to the opposite algebras.
The corresponding arrows in the opposite quiver are
$a_i^{op}:v_{f_i}\to v_{e_i}$, and our hypothesis gives
\[
\mr{Hom}_{\mr{mod}\text{-}\Lambda^{op}}
(f_i\Lambda^{op},e_j\Lambda^{op})
\simeq f_i\Lambda e_j=0
\qquad(1\leq i,j\leq s).
\]
Moreover, these arrows do not occur in the corresponding
minimal generators of the opposite relation ideal.
Thus $\Lambda^{op}\to\Gamma^{op}$ satisfies the
arrow-removal hypotheses. Applying \cite[Proposition~4.6]{GPS2021} to this arrow removal
and identifying left modules with right modules over the
opposite algebras, we obtain a cleft extension
\[
(\Gamma\text{-}\mr{mod},
 \Lambda\text{-}\mr{mod},i,\mathsf{e},l).
\]
Furthermore, $l$ and $r$ are exact, and $\mathsf{e}$ preserves
projective modules.
By \cite[Corollary~3.3]{EPS}, applied in the same way,
$\mathsf{e}$ is a $1$-eventually homological isomorphism.
The following theorem now follows from
Proposition \ref{P-cleft-bound}(b).
\begin{Th}\label{T-arrow-bound}
For an arrow removal  $\Lambda\to\Gamma$ as above,
\begin{align*}
\mr{gb}(\Lambda)&\leq\max\{1,\mr{gb}(\Gamma)\},&
\mr{gb}(\Gamma)&\leq\max\{1,\mr{gb}(\Lambda)\},\\
\mr{fab}(\Lambda)&\leq\max\{1,\mr{fab}(\Gamma)\},&
\mr{fab}(\Gamma)&\leq\max\{1,\mr{fab}(\Lambda)\}.
\end{align*}
\end{Th}

\begin{Rem}\label{R-cleft-bound}
The Uniform Auslander Condition, as well as finiteness of the finitistic Auslander bound, is invariant under arrow removal.
If either member of one of the two pairs of invariants is greater than $1$, then the two members of that pair are equal.
\end{Rem}
The following example illustrates the role of the constant $1$
in Theorem \ref{T-arrow-bound}. Although arrow removal preserves
finiteness of both Auslander bounds, it need not preserve their
values when they are at most $1$. In this example, both bounds
decrease from $1$ to $0$ after removing a single arrow, while
both algebras have infinite global dimension.
\begin{Exm}\label{E-arrow-bound}
Let $Q$ be the quiver
\[
\xymatrix{
1 \ar@(ul,dl)[]_{x} \ar[r]^a & 2
}
\]
and put $\Lambda=\Bbbk Q/\langle x^2\rangle$.
Paths are composed from left to right.
Removing the arrow $a$ gives
\[
\Gamma=\Lambda/\Lambda\overline a\Lambda
\simeq D\times\Bbbk,
\qquad D=\Bbbk[x]/(x^2).
\]
Then
\[
\mr{gb}(\Gamma)=\mr{fab}(\Gamma)=0,
\qquad
\mr{gb}(\Lambda)=\mr{fab}(\Lambda)=1.
\]
Thus the constant $1$ in Theorem \ref{T-arrow-bound}
cannot in general be omitted.
\end{Exm}

\begin{proof}
The arrow $a$ does not occur in the minimal generating
set $\{x^2\}$ of the relation ideal. Moreover,
$
\mr{Hom}_\Lambda(\Lambda e_2,\Lambda e_1)
\simeq e_2\Lambda e_1=0,
$
since there is no path from vertex $2$ to vertex $1$.
Hence $\Lambda\to\Gamma$ is an arrow removal satisfying
the stated hypotheses.

We first calculate the bounds of $\Gamma$. By Example \ref{E-T1-nontriangular},  $ \mr{gb}(D)=\mr{fab}(D)=0.$
Since Ext over a product algebra is computed
componentwise and $\Bbbk$ is semisimple, we obtain
\[
\mr{gb}(\Gamma)=\mr{fab}(\Gamma)=0.
\]
Theorem \ref{T-arrow-bound} therefore yields
$
\mr{gb}(\Lambda)\leq1,$ $\mr{fab}(\Lambda)\leq1.
$ To prove equality, put $P_i=\Lambda e_i$ and
$S_i=P_i/\mr{rad}P_i$.
Right multiplication by $a$ gives an exact sequence
\[
0\longrightarrow P_1
\xrightarrow{\cdot a}P_2
\longrightarrow S_2\longrightarrow0.
\]
Applying $\mr{Hom}_\Lambda(-,S_1)$ and using
$\mr{Hom}_\Lambda(P_2,S_1)=0$ and
$\mr{Hom}_\Lambda(P_1,S_1)\simeq\Bbbk$, we obtain
\[
\mr{Ext}_\Lambda^1(S_2,S_1)\simeq\Bbbk.
\]
Consequently,
\[
\mr{pd}_\Lambda S_2=b_\Lambda(S_2)=1.
\]
Together with the upper bounds, this proves
\[
\mr{gb}(\Lambda)=\mr{fab}(\Lambda)=1.
\]

Finally, both algebras have infinite global dimension. Thus the example also illustrates the arrow-removal
result beyond finite global dimension.
\end{proof}
\section{An arbitrarily large gap between $\mr{fab}(\Lambda)$ and $\findim(\Lambda)$}
\label{S-radical-square-zero}
In this section, we give a family with zero finitistic dimension and arbitrarily large finite
Auslander bounds. 

Fix $n\geq 1$.  Let $Q_n$ be the quiver
\[
\xymatrix{n & n-1 \ar[l]_{\alpha_n} & \cdots \ar[l] & 1 \ar[l] & 0 \ar[l]_{\alpha_1} \ar@(ur,dr)[]^{\varepsilon}}
\]
Let $J$ be the arrow ideal of $kQ_n$ and set $ \Lambda_n=kQ_n/J^2.$  
For each vertex $i$, let $S_i$ be the corresponding simple left module
and $P_i=\Lambda_n e_i$ its projective cover.  Since $J^2=0$, the radicals
of the indecomposable projectives are semisimple, with
\begin{equation}\label{eq:rad}
 \mr{rad} P_0\cong S_0,\qquad
 \mr{rad} P_i\cong S_{i-1}\quad(1\le i\le n).
\end{equation}
In particular,
\begin{equation}\label{eq:syzygies}
 \Omega S_0\cong S_0,\qquad
 \Omega S_i\cong S_{i-1}\quad(1\le i\le n).
\end{equation}
All simple left $\Lambda_n$-modules consequently have infinite
projective dimension. Clearly, $ \findim\Lambda_n=0.$ In the following, we prove that
 \[\mr{fab}(\Lambda_n)=\mr{gb}(\Lambda_n)=n.\]

\begin{Lem}\label{lem:simple-bounds}
For every $0\le i\le n$, $b_{\Lambda_n}(S_i)=i.$
\end{Lem}
\begin{proof}
Since $\Omega S_0\cong S_0$, dimension shifting gives
\[
 \Ext^{j+1}_{\Lambda_n}(S_0,N)
 \cong\Ext^j_{\Lambda_n}(S_0,N)
 \qquad(j\ge1)
\]
for every finitely generated $N$.  If these groups eventually vanish,
they all vanish in positive degrees.  Hence $b_{\Lambda_n}(S_0)=0$.

For $i\ge1$, $\Omega^i S_i\cong S_0$ by
\eqref{eq:syzygies}.  Consequently, if
$\Ext^j_{\Lambda_n}(S_i,N)=0$ for all sufficiently large $j$,
then dimension shifting and the preceding paragraph show $ \Ext^j_{\Lambda_n}(S_i,N)=0\qquad(j>i).$ 
It follows that $b_{\Lambda_n}(S_i)\le i$.

To prove the reverse inequality, take $N=P_0$.
The exact sequence
\[
 0\longrightarrow S_0\longrightarrow P_0
 \longrightarrow S_0\longrightarrow0
\]
shows that $\Ext^1_{\Lambda_n}(S_0,P_0)=0$ since 
 the restriction map $\Hom_{\Lambda_n}(P_0,P_0)\longrightarrow
\Hom_{\Lambda_n}(S_0,P_0)$ 
is surjective.  Periodicity of $S_0$ then yields
\begin{equation}\label{eq:S0-P0}
 \Ext^j_{\Lambda_n}(S_0,P_0)=0\qquad(j\ge1).
\end{equation}
For $S_1$, the projective cover sequence is
\[
 0\longrightarrow S_0\longrightarrow P_1
 \longrightarrow S_1\longrightarrow0.
\]
Now
$\Hom_{\Lambda_n}(P_1,P_0)=0$, whereas
$\Hom_{\Lambda_n}(S_0,P_0)\neq0$.  Therefore
\[
 \Ext^1_{\Lambda_n}(S_1,P_0)\neq0,
 \qquad
 \Ext^{>1}_{\Lambda_n}(S_1,P_0)=0,
\]
where the second assertion follows from \eqref{eq:S0-P0}.
For arbitrary $i\ge1$, since $\Omega^{i-1}S_i\cong S_1$,
dimension shifting gives
\[
 \Ext^i_{\Lambda_n}(S_i,P_0)
 \cong\Ext^1_{\Lambda_n}(S_1,P_0)\neq0,
 \qquad
 \Ext^{>i}_{\Lambda_n}(S_i,P_0)=0.
\]
Thus the last nonzero extension group for $(S_i,P_0)$ occurs in degree $i$, proving
$b_{\Lambda_n}(S_i)\ge i$ and hence equality.
\end{proof}

\begin{Pro}\label{thm:main}
$\mr{fab}(\Lambda_n)=\mr{gb}(\Lambda_n)=n.$
\end{Pro}
\begin{proof}
Let $M\in\Lambda_n\text{-}\mr{mod}$ be nonprojective.  Then $ \Omega M\cong\bigoplus_{i=0}^{n-1}S_i^{a_i}$ 
for suitable nonnegative integers $a_i$. If $\Ext^{\gg0}_{\Lambda_n}(M,N)=0$, then
$\Ext^{\gg0}_{\Lambda_n}(S_i,N)=0$ for every $i$ with $a_i>0$.
Lemma~\ref{lem:simple-bounds} implies
$\Ext^{>i}_{\Lambda_n}(S_i,N)=0$ for such $i$.
Dimension shifting therefore gives
\[
 \Ext^{j}_{\Lambda_n}(M,N)
 \cong\Ext^{j-1}_{\Lambda_n}(\Omega M,N)=0
 \qquad(j>n).
\]
Hence $b_{\Lambda_n}(M)\le n$.  Thus $\mr{gb}(\Lambda_n)\le n$.
On the other hand, Lemma~\ref{lem:simple-bounds} gives
$b_{\Lambda_n}(S_n)=n$.  Therefore
\[
 n\le\mr{fab}(\Lambda_n)\le\mr{gb}(\Lambda_n)\le n,
\]
which proves the result.
\end{proof}

\section*{Acknowledgements}
Kaili Wu is supported by  the  Natural Science Foundation of Jiangsu Province of China (Grant No. BK20240657) and the National Natural Science
Foundation of China (Grant No. 12526627).

 Jiaqun Wei is supported by the National Natural Science Foundation of China (Grant Nos. 12571042, 12271248) and the Natural Science
Foundation of Zhejiang Province (Grant No. LZ25A010002).\\
\section*{Conflict of Interest}  The authors declare that they have no conflict of interest.

\end{document}